\documentclass[11pt]{article}
\begin{document}


\newtheorem{thm}{Theorem}[section]


\newtheorem{lem}[thm]{Lemma}
\newtheorem{cor}[thm]{Corollary}
\newtheorem{ex}[thm]{Example}
\newtheorem{prop}[thm]{Proposition}
\newtheorem{remark}[thm]{Remark}
\newtheorem{coun}[thm]{Counterexample}
\newtheorem{defn}[thm]{Definition}
\newtheorem{conj}{Conjecture}
\newtheorem{problem}{Problem}
\newcommand\ack{\section*{Acknowledgement.}}

\newcommand{\etal}{{\it et al. }}


\newcommand{\bbP}{{\rm I\hspace{-0.8mm}P}}
\newcommand{\bbE}{{\rm I\hspace{-0.8mm}E}}
\newcommand{\bbF}{{\rm I\hspace{-0.8mm}F}}
\newcommand{\bbI}{{\rm I\hspace{-0.8mm}I}}
\newcommand{\bbR}{{\rm I\hspace{-0.8mm}R}}
\newcommand{\bbRp}{{\rm I\hspace{-0.8mm}R}_+}
\newcommand{\bbN}{{\rm I\hspace{-0.8mm}N}}
\newcommand{\bbC}{{\rm C\hspace{-2.2mm}|\hspace{1.2mm}}}
\newcommand{\bbD}{{\rm I\hspace{-0.8mm}D}}
\newcommand{\bbQ}{\bf Q}
\newcommand{\bbZ}{{\rm \rlap Z\kern 2.2pt Z}}
\newcommand{\bbK}{{\rm I\hspace{-0.8mm}K}}

\newcommand{\matP}{{\bbP}}
\newcommand{\mattildeP}{\tilde{\bbP}}
\newcommand{\matPN}[1]{{\bbP}_{#1}}
\newcommand{\matPP}[1]{{\bbP}_{#1}^0}
\newcommand{\matE}{{\bbE}}
\newcommand{\mattildeE}{\tilde{\bbE}}
\newcommand{\matEP}[1]{{\bbE}_{#1}^0}
\newcommand{\matF}{{\bbF}}
\newcommand{\matR}{{\bbR}}
\newcommand{\matRp}{{\bbRp}}
\newcommand{\matN}{{\bbN}}
\newcommand{\matZ}{{\bbZ}}
\newcommand{\matI}{{\bbI}}
\newcommand{\matK}{{\bbK}}
\newcommand{\matQ}{{\bbQ}}
\newcommand{\matC}{{\bbC}}
\newcommand{\matD}{{\bbD}}

\newcommand{\calL}{{\cal L}}
\newcommand{\calM}{{\cal M}}
\newcommand{\calN}{{\cal N}}
\newcommand{\calF}{{\cal F}}
\newcommand{\calG}{{\cal G}}
\newcommand{\calD}{{\cal D}}
\newcommand{\calB}{{\cal B}}
\newcommand{\calH}{{\cal H}}
\newcommand{\calI}{{\cal I}}
\newcommand{\calP}{{\cal P}}
\newcommand{\calQ}{{\cal Q}}
\newcommand{\calS}{{\cal S}}
\newcommand{\calT}{{\cal T}}
\newcommand{\calC}{{\cal C}}
\newcommand{\calK}{{\cal K}}
\newcommand{\calX}{{\cal X}}
\newcommand{\cals}{{\cal S}}
\newcommand{\calE}{{\cal E}}

\newcommand{\koniecmat}{\,}

\newcommand{\eqd}{\ =_{\rm d}\ }
\newcommand{\toto}{\leftrightarrow}
\newcommand{\eqdistr}{\stackrel{\rm d}{=}}
\newcommand{\as}{\stackrel{\rm a.s.}{=}}
\newcommand{\assubs}{{\rm a.s.}}
\newcommand{\convdistr}{\stackrel{\rm d}{\rightarrow}}
\newcommand{\convweak}{{\Rightarrow}}
\newcommand{\convas}{\stackrel{\rm a.s.}{\rightarrow}}
\newcommand{\convfidi}{\stackrel{\rm fidi}{\rightarrow}}
\newcommand{\convprob}{\stackrel{p}{\rightarrow}}
\newcommand{\deff}{\stackrel{\rm def}{=}}
\newcommand{\bis}{{'}{'}}
\newcommand{\Cov}{{\rm Cov}}
\newcommand{\Var}{{\rm Var}}
\newcommand{\Exp}{{\rm E}}

\newcommand{\nd}{n^{\delta}}
\newcommand{\koniec}{\newline\vspace{3mm}\hfill $\odot$}


\title{Weak convergence of Vervaat and Vervaat Error processes of long-range dependent sequences  \protect}
\author{Mikl\'{o}s Cs\"{o}rg\H{o}
\thanks{School of Mathematics and Statistics, Carleton University, 1125 Colonel By Drive, Ottawa, Ontario, K1S 5B6 Canada, email:
mcsorgo@math.carleton.ca}
\and Rafa{\l} Kulik
\thanks{School of Mathematics and Statistics, University of Sydney, NSW 2006, Australia, email: rkuli@maths.usyd.edu.au and Mathematical Institute, Wroc{\l}aw University, Pl. Grunwaldzki 2/4, 50-384 Wroc{\l}aw, Poland.
\newline Research supported in part by NSERC Canada Discovery Grants of
Mikl\'{o}s Cs\"{o}rg\H{o}, Donald Dawson and Barbara Szyszkowicz at
Carleton University}}
\maketitle


\begin{abstract}
Following Cs\"{o}rg\H{o}, Szyszkowicz and Wang (Ann. Statist. {\bf
34}, (2006), 1013--1044) we consider a long range dependent linear
sequence. We prove weak convergence of the uniform Vervaat and the
uniform Vervaat error processes, extending their results to
distributions with unbounded support and removing normality
assumption.
\end{abstract}



\section{Introduction}
Let $\{\epsilon_i,i\in {\matZ}\}$ be a centered sequence of i.i.d.
random variables with finite variance. Consider the class of
stationary linear processes
\begin{equation}\label{model}
X_i=\sum_{k=0}^{\infty}c_k\epsilon_{i-k},\ \ \ i\ge 1 .
\end{equation}
We assume that the sequence $c_k$, $k\ge 0$, is regularly varying
with index $-\beta$, $\beta\in (1/2,1)$ (written as $c_k\in
RV_{-\beta}$). This means that $c_k\sim k^{-\beta}L_0(k)$ as
$k\to\infty$, where $L_0$ is slowly varying at infinity (see e.g.
\cite[Sections 1.4, 1.5]{BinghamGoldieTeugesl1987} for the
definition of slowly varying functions). We shall refer to all such
models as long range dependent (LRD) linear processes. In
particular, by the Karamata Theorem, the covariances $\rho_k:=\Exp
X_0X_k$ decay at the hyperbolic rate, $\rho_k=L(k)k^{-(2\beta-1)}$,
where $\lim_{k\to\infty}L(k)/L_0^2(k)=B(2\beta-1,1-\beta)$ and
$B(\cdot,\cdot)$ is the beta-function. Consequently, since
$-(2\beta-1)>-1$, the covariances are not summable.

Assume that $X_1$ has a continuous distribution function $F$ and the
density $f$, which is assumed to be positive almost everywhere. For
$y\in (0,1)$ define $Q(y)=\inf\{x:F(x)\ge y\}=\inf\{x:F(x)= y\}$,
the corresponding (continuous) quantile function. Given the ordered
sample $X_{1:n}\le\cdots\le X_{n:n}$ of $X_1,\ldots,X_n$, let
$F_n(x)=n^{-1}\sum_{i=1}^n1_{\{X_i\le x\}}$ be the empirical
distribution function and $Q_n(\cdot)$ be the corresponding
left-continuous sample quantile function. Define $U_i=F(X_i)$ and
$E_n(x)=n^{-1}\sum_{i=1}^n1_{\{U_i\le x\}}$, the associated uniform
empirical distribution. Denote by $U_n(\cdot)$ the corresponding
uniform sample quantile function.

Let $r$ be an integer and let
$$Y_{n,r}=\sum_{i=1}^n\sum_{1\le j_1<\cdots<
j_r<\infty}\prod_{s=1}^rc_{j_s}\epsilon_{i-j_s},\qquad n\ge 1,$$ so
that $Y_{n,0}=n$, and $Y_{n,1}=\sum_{i=1}^nX_i$. If $1\le
p<(2\beta-1)^{-1}$, then (cf. \cite{HoHsing1996})
\begin{equation}\label{eq-varance-behaviour}
\sigma_{n,p}^2:=\Var (Y_{n,p})\sim const. n^{2-p(2\beta-1)}L^{2p}(n).
\end{equation}
In particular
$$
\sigma_{n,1}^2\sim \frac{c_{\beta}}{(1-\beta)(3-2\beta)}n^{3-2\beta}L^2(n)=:n^{3-2\beta}L_0^2(n).
$$

Define now the general empirical, the uniform empirical, the general
quantile and the uniform quantile processes respectively as follows:
\begin{equation}\label{eq.general.empirical}
\beta_n(x)=\sigma_{n,1}^{-1}n(F_{n}(x)-F(x)),\qquad x\in {\matR},
\end{equation}
\begin{equation}\label{eq.uniform.empirical}
\alpha_n(y)=\sigma_{n,1}^{-1}n(E_{n}(y)-y),\qquad y\in [0,1],
\end{equation}
\begin{equation}\label{eq.general.quantile}
q_n(y)=\sigma_{n,1}^{-1}n(Q(y)-Q_{n}(y)),\qquad y\in (0,1),
\end{equation}
\begin{equation}\label{eq.uniform.quantile}
u_n(y)=\sigma_{n,1}^{-1}n(y-U_{n}(y)),\qquad y\in [0,1].
\end{equation}
\noindent Let
\begin{equation}\label{eq-definition-BahadurKiefer-uniform}
\tilde{R}_n(y)=\alpha_n(y)-u_n(y),\qquad y\in [0,1],
\end{equation}
be the uniform Bahadur-Kiefer process. This process was introduced by Kiefer in \cite{Kiefer1970}, though not explicitly,
in order to study the behavior of quantile processes via that of empirical, as initiated by Bahadur \cite{Bahadur1966} for $y\in (0,1)$ fixed.\\

Let
$$\tilde
V_n(t)=2\sigma_{n,1}^{-1}n\int_0^t\tilde R_n(y)dy,\qquad t\in
[0,1],$$ be the uniform Vervaat process and
$$
\tilde W_n(t)=2\sigma_{n,1}^{-1}n\int_0^t\tilde
R_n(y)dy-\alpha_n^2(t), \qquad t\in [0,1],
$$
be the uniform Vervaat error process as in \cite{CSW2006}.\\

Assume for a while that $\{\eta_n\}_{n\ge 1}$ is a stationary and
standardized (i.e., zero-mean and unit variance) long-range
dependent Gaussian sequence with a covariance structure
$$
\gamma (k):=\Exp (\eta_1\eta_{k+1})=k^{-D}\tilde{L}(k), \qquad 0<D<1,
$$
where $\tilde{L}$ is slowly varying at infinity. Let $G$ be an
arbitrary real-valued measurable function and define
$Y_n=G(\eta_n)$, $n\ge 1$. Let $F_Y$ be the continuous distribution
function of $Y_1$ and $Q_Y(\cdot)$ the corresponding continuous
quantile function. Define $V_i=F_Y(Y_i)$. As in Dehling and Taqqu
\cite{DehlingTaqqu1989}, expand $1_{\{X_n\le x\}}-F(x)$ as,
$$
1_{\{Y_n\le x\}}-F_Y(x)=\sum_{l=\tau_x}^{\infty}c_l(x)H_l(\eta_n)/l!,
$$
where $$H_l(x)=(-1)^l\exp(x^2/2)\frac{d^l}{dx^l}\exp(-x^2/2)$$ is
the $l$th Hermite polynomial,
$$c_l(x)=\Exp \left[\left(1_{\{G(\eta_1)\le
x\}}-F_Y(x)\right)H_l(\eta_1)\right],$$ and for any $x\in {\matR}$,
$\tau_x$ (the Hermite rank) is the index of the first non-zero
coefficient of the expansion. The uniform version is obtained as
$$
1_{\{V_n\le y\}}-y=\sum_{l=\tau_y}^{\infty}J_l(y)H_l(\eta_n)/l!,
$$
where now $J_{l}(y)=c_l(Q_Y(y))$ for any $y\in (0,1)$.

Let $\tilde{\sigma}_{n,\tau}^2=n^{2-\tau D}\tilde {L}^{\tau}(n)$.
Replace the constants $\sigma_{n,1}$ with $\tilde\sigma_{n,\tau}$ in
the definitions of $\tilde R_n(\cdot)$, $\tilde V_n(\cdot)$ and
$\tilde W_n(\cdot)$. In \cite{CSW2006} Cs\"{o}rg\H{o}, Szyszkowicz
and Wang (CsSzW) proved that the uniform Bahadur-Kiefer process
$\tilde R_n(\cdot)$ converges weakly in $D([0,1])$. This phenomenon
is exclusive for long range dependent sequences, since in the i.i.d.
case the (uniform) Bahadur-Kiefer process cannot converge weakly.
However, as it was first shown by Vervaat \cite{Vervaat1972}, in the
i.i.d case the uniform Vervaat process does converge weakly.
Obviously, in the LRD case, weak convergence of the uniform Vervaat
process is implied by that of $\tilde R_n(\cdot)$, namely (see
\cite[Theorem 3.1]{CSW2006}):
\begin{equation}\label{weak-conv-uniform-verwaat-CSW}
\tilde V_n(t)\convweak \frac{2}{(2-\tau D)(1-\tau
D)}J_{\tau}^2(t)Z_{\tau}^2 ,\qquad n\to\infty,
\end{equation}
where $\convweak$ denotes weak convergence in $D([0,1])$ equipped
with the sup-norm, and $Z_{\tau}$ is a random variable defined by an
appropriate integral with respect to Brownian motion (see
\cite{DehlingTaqqu1989}). In particular, if $\tau=1$, then $Z_1$ is
standard normal. Further, CsSzW \cite{CSW2006} observed that,
similarly to the i.i.d case, the limiting process associated with
$\tilde V_n(\cdot)$ agrees with that of $\alpha_n^2(\cdot)$.
Therefore, it makes sense to consider the uniform Vervaat error
process $\tilde W_n(\cdot)$. They showed that this process converges
weakly as well, via concluding
\begin{equation}\label{weak-conv-uniform-verwaat-error-CSW}
n\sigma_{n,1}^{-1}\tilde W_n(t)\convweak \frac{2^{5/2}}{(2-\tau
D)^{3/2}(1-\tau D)^{3/2}}J_{\tau}^2(t)J_{\tau}^{'}(t)Z_{\tau}^3,
\qquad n\to\infty .
\end{equation}
This property is also exclusive for the LRD case. We refer to
\cite{CaskiCsorgoFoldesShiZitikise}, \cite{CsorgoZitikis1999},
\cite{Zitikis1998}, \cite{CsorgoSzyszkowicz1998} as well as the
Introduction in \cite{CSW2006} for motivations, probabilistic
properties and applications of Bahadur-Kiefer, Vervaat and Vervaat
error processes.

We note in passing that, though the results in CsSzW \cite{CSW2006}
for the uniform Bahadur-Kiefer process and, consequently, for the
uniform Vervaat and Vervaat error processes, are true, their proofs
are invalid, unless $F_Y$, the distribution of the subordinated
random variable $G(\eta_1)$, is assumed to have finite support.
Moreover, even then, the limiting process in
(\ref{weak-conv-uniform-verwaat-error-CSW}) should be corrected via
multiplying it by $\frac{1}{2}$, see \cite{CSW2006-corr}. \\

In case of the Bahadur-Kiefer process, the problem of an infinite
support was solved in \cite{CsorgoKulik2006} in a more general
setting in the case of LRD linear sequences by using weighted
approximations. However, in general, this is still not suitable for
establishing the weak convergence of the Vervaat process $\tilde
V_n(\cdot)$, unless some specific conditions are imposed on the
model. The reason for the problems arising in \cite{CSW2006}, and
faced up to in \cite{CsorgoKulik2006}, is that, unlike in the i.i.d.
case, the uniform quantile process contains information about the
quantile
function associated with the random variables $X_n$.\\

Therefore, coming back to LRD linear sequences, the aim of this
paper is to present an appropriate approximation result for the
uniform Bahadur-Kiefer process, which will be suitable to treat the
uniform Vervaat process to obtain
(\ref{weak-conv-uniform-verwaat-CSW}), when $F$ is assumed to have
infinite support . Further, we will obtain the correct version of
the weak convergence of the uniform Vervaat error process. The
approach is via weighted approximation of the Bahadur-Kiefer process
like in \cite{CsorgoKulik2006}. Thus, first we get the correct
limiting behaviour of the Vervaat error process, second, we remove
assumptions on bounded support of $F$, third, we remove the
normality assumption on $\epsilon_i$. This approach in fact requires
very precise knowledge on the behavior of the
density-quantile function $f(Q(y))$.\\

However, we do not extend the results in \cite{CSW2006} in full
generality, since we do not consider subordinated LRD sequences
$Y_i=G(X_i)$, $i\ge 1$, where $G$ is a measurable function. If $G$
has a power rank 1 (see e.g. \cite{Hsing}), then in expense of some
additional technicalities, the results will be similar as for
non-subordinated case. However, if the power rank is greater than 1,
the scaling
factors and the limiting processes will be different. \\

To state our results, Let $F_{\epsilon}$ be the distribution
function of the centered i.i.d. sequence $\{\epsilon_i,i\in
{\matZ}\}$ with finite $4$th moment. Assume that for a given integer
$p$, the derivatives $F^{(1)}_{\epsilon}, \ldots,
F^{(p+3)}_{\epsilon}$ of $F_{\epsilon}$ are bounded and integrable.
Note that these properties are inherited by the distribution $F$ as
well (cf. \cite{Wu2003}). These conditions will be assumed throughout the paper with $p=2$.\\

We shall need the following conditions on $fQ(\cdot)=f(Q(\cdot))$ and $f^{'}Q(\cdot)=f^{'}(Q(\cdot))$:
\begin{itemize}
\item[{\rm (A)}] $\sup_{y\in (0,1)}|gQ(y)|/(y(1-y))^{1-\mu}=O(1)$ for
some $1/2>\mu>0$ and $g=f,f^{'}$;
\item[{\rm (B)}] $\sup_{y\in (0,1)}|(gQ(y))^{'}|(y(1-y))^{\mu}=O(1)$ for any
$\mu>0$ and $g=f,f^{'}$. Note that
$(f(Q(y)))^{'}f(Q(y))=f^{'}(Q(y))$;
\item[{\rm (C)}] $\sup_{y\in (0,1)}|(gQ(y))^{''}|(y(1-y))^{1+\mu}=O(1)$ for any
$\mu>0$ and $g=f,f^{'}$.
\end{itemize}
We shall prove the following results.
\begin{thm}\label{thm-conv-uniform-BK}
Assume that conditions {\rm(A)-(C)} are fulfilled and $\beta<3/4$.
Then, as $n\to\infty$,
$$
\sup_{y\in [\delta_n,1-\delta_n]}\left|n\sigma_{n,1}^{-1}\tilde
R_n(y)-\sigma_{n,1}^{-2}f^{'}(Q(y))\left(\sum_{i=1}^nX_i\right)^2\right|=o_{\assubs}(1),
$$
where $\delta_n=Cn^{-(2\beta-1)}L_0^{2}(n)(\log\log n)$.
\end{thm}
\begin{cor}\label{thm-conv-uniform-BK-corollary}
Under the conditions of {\rm Theorem}
{\rm\ref{thm-conv-uniform-BK}}, as $n\to\infty$,
$$
n\sigma_{n,1}^{-1}\tilde R_n(y)1_{\{y\in
[\delta_n,1-\delta_n]\}}\convweak f^{'}(Q(y))Z_{1}^2 .
$$
\end{cor}
\begin{thm}\label{thm-conv-uniform-verwaat}
Under the conditions of {\rm Theorem}
{\rm\ref{thm-conv-uniform-BK}}, as $n\to\infty$,
$$\tilde V_n(t)\convweak f^2Q(t)Z_{1}^2.
$$
\end{thm}
\begin{thm}\label{thm-conv-uniform-verwaat.error}
Under the conditions of {\rm Theorem}
{\rm\ref{thm-conv-uniform-BK}}, as $n\to\infty$,
\begin{equation}\label{eq-VE}
\sigma_{n,1}^{-1}n\tilde W_n(t)\convweak \frac{1}{((3-
2\beta)(1-\beta))^{3/2}}f^2Q(t)(fQ)^{'}(t)Z_1^3.
\end{equation}
\end{thm}
\begin{remark}{\rm
A few words on the conditions (A)-(C). Assume that $F=\Phi$ (the
standard normal distribution). It follows from \cite{Parzen1979}
that (A) is fulfilled. Further,
$(\phi(\Phi^{-1}(y)))^{'}=-\Phi^{-1}(y)$ is unbounded (this is
actually the reason, why the proofs in \cite{CSW2006} do not work),
but (B) holds. Furthermore,
$(\phi(\Phi^{-1}(y)))^{''}=-\frac{1}{\phi(\Phi^{-1}(y))}$, and it
follows from \cite{Parzen1979} that (C) is fulfilled.

Furthermore, one can check that the conditions (A)-(B) are fulfilled
for distributions with exponential or Pareto tails. To be more
specific, let $f(x)=const. |x|^{-\alpha}$, $x>\delta>0$, $\alpha>2$.
Also, in $[-\delta,\delta]$, $f$ is interpolated smoothly to assure
existence of its derivatives - note that most important issue in
(A)-(C) is the tail behaviour of the density. Then, for $x>\delta$,
$F(x)=1-c x^{-(\alpha-1)}$, $c\in (0,\infty)$, and
$Q(y)=c^{1/(\alpha-1)}(1-y)^{-1/(\alpha-1)}$ for $y>\delta_0>1/2$.
Consequently, $f(Q(y))/(1-y)^{1-\mu}=const.
(1-y)^{\mu+\frac{1}{\alpha-1}}$ and $\sup_{y>\delta_0}
f(Q(y))/(1-y)^{1-\mu}=O(1)$. Also, $\sup_{y>\delta_0}
f^{'}(Q(y))/(1-y)^{1-\mu}=O(1)$. The similar consideration applies
to the left tail. Consequently, the condition (A) is fulfilled.
Conditions (B) and (C) can be verified in a similar way.

More generally, if $f_{\epsilon}(x)=|x|^{-\alpha}L_1(x)$, $L_1$
being slowly varying at infinity, then
$\lim_{x\to\infty}P(X_1>x)/P(|\epsilon_1|>x)=const.\in (0,\infty)$
(see e.g. \cite{MikoschSam}) and by the Karamata Theorem,
$\lim_{x\to\infty} f(x)/(x^{-\alpha}L_1(x))=const. \in (0,\infty)$.
Recalling that (A)-(C) are essentially the conditions on the
asymptotic tail behaviour, we conclude that (A)-(C) hold.
 }
\end{remark}
\begin{remark}{\rm
In Theorem \ref{thm-conv-uniform-BK} we are not able to obtain the
a.s. approximation on $(0,1)$. From this theorem, weak convergence
of $\tilde R_n(y)1_{\{y\in [\delta_n,1-\delta_n]\}}$ follows, as in
Corollary \ref{thm-conv-uniform-BK-corollary}. We are not able to
obtain weak convergence on $(0,1)$ either. However, this was not our
concern in this paper. It can be done via weight functions (see
\cite{CsorgoKulik2006} for more details). Nevertheless, this
convergence is good enough to obtain weak convergence of both the
uniform Vervaat and the uniform Vervaat error processes. The weak
convergence limit in Theorem \ref{thm-conv-uniform-verwaat.error}
differs from that of Proposition 3.2 in \cite{CSW2006} by the
already mentioned factor of $\frac{1}{2}$. To see this, assume that
$\Exp(\epsilon_1^2)=1$ and note that parametrization of the Gaussian
and the linear model yields $\tilde L(n)=c_{\beta}L^2(n)$,
$D=2\beta-1$. Plugging this into (\ref{eq-VE}) we see, that the
result (\ref{weak-conv-uniform-verwaat-error-CSW}) should be
corrected by replacing $2^{5/2}$ with $2^{3/2}$.

The problem in the proof
of Proposition 3.2 in \cite{CSW2006} comes from an inappropriate use
of their Proposition 2.5.}
\end{remark}

In what follows $C$ will denote a generic constant which may be
different at each time it appears. Further, $\ell(n)$ is a slowly varying function at infinity, possibly different at each time it appears.\\

\section{Proofs}
Recall that
$$
\delta_n=n^{-(2\beta-1)}L_0^{2}(n)(\log\log n)
$$
and let
$$
a_n=n^{-(\beta-1/2)}L_0(n)(\log\log n)^{1/2},
$$
$$
d_{n,p}=\left\{\begin{array}{ll} n^{-(1-\beta)}L_0^{-1}(n)(\log n)^{5/2}(\log\log n)^{3/4}, & (p+1)(2\beta-1)>1\\
n^{-p(\beta-\frac{1}{2})}L_0^{p}(n)(\log n)^{1/2}(\log\log n)^{3/4},
& (p+1)(2\beta-1)<1
\end{array}\right. ,
$$
Note that
$d_{n,2}=o(a_n)$ if $\beta<3/4$ and $\sigma_{n,1}^{-1}=o(d_{n,2})$.

\subsection{Preliminary results}
We recall the following law of the iterated logarithm for partial
sums $\sum_{i=1}^nX_i$ (see, e.g., \cite{WangLinGulati}):
\begin{equation}\label{step.5}
\limsup_{n\to\infty}\sigma_{n,1}^{-1}(\log\log
n)^{-1/2}\left|\sum_{i=1}^nX_i\right|\as c(\beta,1),
\end{equation}
where where
$c^2(\beta,p)=\left(\int_0^{\infty}x^{-\beta}(1+x)^{-\beta}dx\right)(1-\beta)^{-1}(3-2\beta)^{-1}$.
Also, if $1\le p<(2\beta-1)^{-1}$, then
\begin{equation}\label{eq-conv-partial.sums}
Y_{n,p}=O_P(\sigma_{n,p}).
\end{equation}
\begin{lem}\label{lem.strong.bound}
Let $p\ge 1$ be an arbitrary integer such that $p<(2\beta-1)^{-1}$.
Then, as $n\to\infty$,
\begin{equation}\label{conjecture}
Y_{n,p}=O_{\assubs}(\sigma_{n,p}(\log n)^{1/2}\log\log n).
\end{equation}
\end{lem}
{\it Proof.} Let $B_n^2=\sigma_{n,p}^2\log n (\log\log n)^2$. By
(\ref{eq-varance-behaviour}), \cite[Lemma 4]{Wu2005} and Karamata's
Theorem we have for $2^{d-1}<n\le 2^d$,
\begin{eqnarray*}
\lefteqn{\hspace*{-2cm} \left\|\max_{k\le
n}\frac{Y_{k,p}}{B_{2^d}}\right\|_2^2\le
\frac{1}{B_{2^d}}\left(\sum_{j=0}^{d}2^{(d-j)/2}\sigma_{2^j,p}\right)^2\le
\frac{1}{B_{2^d}^2}\left(\sum_{j=0}^d2^{j(1-p(2\beta-1))/2}L_0^p(2^j)\right)^2}\\
&\sim & \frac{1}{B_{2^d}^2}2^{2d-dp(2\beta-1)}L_0^{2p}(2^d)\sim
d^{-1}(\log d)^{-2}.
\end{eqnarray*}
Therefore, the result follows by the Borel-Cantelli lemma.
 \koniec

The next result gives the reduction principle for the empirical
processes.
\begin{thm}[\cite{Wu2003}]\label{thm-HoHsing}
Let $p$ be a positive integer. Then, as $n\to\infty$,
$$
\Exp \sup_{x\in {\matR}}\left|\sum_{i=1}^n(1_{\{X_i\le
x\}}-F(x))+\sum_{r=1}^p(-1)^{r-1}F^{(r)}(x)Y_{n,r}\right|^2=O(\Xi_n+n(\log
n)^2),
$$
where
$$
\Xi_n=\left\{\begin{array}{ll} O(n), & (p+1)(2\beta-1)>1\\
O(n^{2-(p+1)(2\beta-1)}L_0^{2(p+1)}(n)), & (p+1)(2\beta-1)<1
\end{array}\right. .
$$
\end{thm}

\indent Let $V_{n,p}(x)=\sum_{r=1}^p(-1)^{r-1}F^{(r)}(x)Y_{n,r}$,
$x\in {\matR}$ and $\tilde V_{n,p}(y)=V_{n,p}(Q(y))$, $y\in (0,1)$.
Using Theorem \ref{thm-HoHsing} and the same argument as in the
proof of Lemma \ref{lem.strong.bound}, we obtain
\begin{eqnarray}
\lefteqn{\sup_{x\in {\matR}}|\beta_n(x)+\sigma_{n,1}^{-1}V_{n,p}(x)|=}\label{eq-reduction-principle}\\
& = & \sigma_{n,1}^{-1}\sup_{x\in
{\matR}}\left|\sum_{i=1}^n(1_{\{X_i\le x\}}-F(x))+
V_{n,p}(x)\right|=o_{\assubs}(d_{n,p})\nonumber.
\end{eqnarray}
Consequently, via $\{\alpha_n(y),y\in (0,1)\}=\{\beta_n(Q(y)),y\in
(0,1)\}$,
\begin{equation}\label{approx-unif-empirical}
\sup_{y\in (0,1)}|\alpha_n(y)+\sigma_{n,1}^{-1}\tilde
V_{n,p}(y)|=O_{\assubs}(d_{n,p}).
\end{equation}
We shall use this result with $p=2$. Then, as mentioned before, $d_{n,2}=o(a_n)$ if $\beta<3/4$.
\subsection{Results on the uniform empirical and quantile processes}
We have
$$
\frac{\tilde
V_{n,2}(y)}{(y(1-y))^{1/2}}=\frac{f(Q(y))}{(y(1-y))^{1/2}}\sum_{i=1}^nX_i-\frac{f^{'}(Q(y))}{(y(1-y))^{1/2}}Y_{n,2}.
$$
Write
$$
\frac{f^{'}(Q(y))}{(y(1-y))^{1/2}}=\frac{f^{'}(Q(y))}{f(Q(y))}(y(1-y))^{\mu}
\frac{f(Q(y))}{(y(1-y))^{1/2+\mu}}
$$
with $\mu<1/2$. Using (A), (B) and (\ref{step.5}) we have
$\frac{\tilde V_{n,2}(y)}{(y(1-y))^{1/2}}=O_{\assubs}((\log\log
n)^{1/2})$ uniformly on $(0,1)$.
\begin{lem}\label{lem-weighted.lil-unifom-empirical}
Under the conditions of {\rm Theorem}
{\rm\ref{thm-conv-uniform-BK}}, as $n\to\infty$,
$$
\sup_{y\in
[\delta_n,1-\delta_n]}\frac{|\alpha_n(y)|}{\sqrt{y(1-y)}}=O_{\assubs}((\log\log
n)^{1/2}).
$$
\end{lem}
{\it Proof.} We have
\begin{eqnarray*}
\lefteqn{\sup_{y\in
[\delta_n,1-\delta_n]}\frac{|\alpha_n(y)|}{\sqrt{y(1-y)}}}\\&\le&
\sup_{y\in
[\delta_n,1-\delta_n]}\frac{|\alpha_n(y)+\sigma_{n,1}^{-1}\tilde
V_{n,2}(y)|}{\sqrt{y(1-y)}}+O_{\assubs}((\log\log
n)^{1/2})\\
&= & O_{\assubs}(\delta_n^{-1/2}d_{n,2})+O_{\assubs}((\log\log
n)^{1/2})=O_{\assubs}((\log\log n)^{1/2}),
\end{eqnarray*}
using (\ref{approx-unif-empirical}).\koniec

Using the method of \cite[Theorem 2]{CsorgoRevesz1978}, we obtain
the same result for the uniform quantile process.
\begin{lem}\label{lem-weighted.lil-unifom-quantile}
Under the conditions of {\rm Theorem}
{\rm\ref{thm-conv-uniform-BK}}, with some $C_0\in (0,\infty)$, as $n\to\infty$,
$$
\sup_{y\in
[C_0\delta_n,1-C_0\delta_n]}\frac{|u_n(y)|}{\sqrt{y(1-y)}}=O_{\assubs}((\log\log
n)^{1/2}).
$$
\end{lem}
Next, we study the distance between the empirical and quantile
processes.
\begin{lem}\label{lem-distance-uniform-empirical-quantile}
Under the conditions of {\rm Theorem}
{\rm\ref{thm-conv-uniform-BK}}, as $n\to\infty$,
$$
\sup_{y\in (0,1)}|u_n(y)-\alpha_n(y)|=O_{\assubs}(a_n(\log \log
n)^{1/2}).
$$
\end{lem}
{\it Proof.} Since
\begin{equation}\label{eq-replace}
E_n(U_n(y))=y+O(1/n),
\end{equation}
we obtain from (\ref{approx-unif-empirical}),
\begin{eqnarray*}
\lefteqn{\sup_{y\in [C_0\delta_n,1-C_0\delta_n]}|u_n(y)-\alpha_n(y)|}\\
&\le & \sigma_{n,1}^{-1}\sup_{y\in
[C_0\delta_n,1-C_0\delta_n]}|\tilde
V_{n,2}(y)-\tilde V_{n,2}(U_n(y))|+O_{\assubs}(\sigma_{n,1}^{-1})\\
&\le &
\sigma_{n,1}^{-1}\left|\sum_{i=1}^nX_i\right|\sup_{y\in
[C_0\delta_n,1-C_0\delta_n]}|(fQ)^{'}(\theta)||y-U_n(y)|\\
 && + \sigma_{n,1}^{-1}|Y_{n,2}|\sup_{y\in
[C_0\delta_n,1-C_0\delta_n]}|(f^{'}Q)^{'}(\theta)||y-U_n(y)|+O_{\assubs}(\sigma_{n,1}^{-1}),
\end{eqnarray*}
where $\theta=\theta (y,n)$ is such that $|\theta-y|\le
\sigma_{n,1}n^{-1}|u_n(y)|=O_{\assubs}((y(1-y)\sigma_{n,1}n^{-1}\log\log
n)^{1/2})$ by Lemma \ref{lem-weighted.lil-unifom-quantile}.

Now, via Lemma \ref{lem-weighted.lil-unifom-quantile},
\begin{eqnarray*}
\lefteqn{\sup_{y\in
[\delta_n,1-\delta_n]}|(fQ)^{'}(\theta)||y-U_n(y)|}\\
&=  & \frac{\sigma_{n,1}}{n}\sup_{y\in
[\delta_n,1-\delta_n]}|(fQ)^{'}(\theta)||u_n(y)|\\
&\le &
 \sup_{y\in
[\delta_n,1-\delta_n]}|(fQ)^{'}(\theta)|\sqrt{y(1-y)}O_{\assubs}\left(\frac{\sigma_{n,1}}{n}(\log\log
n)^{1/2}\right)
\end{eqnarray*}
and the bound is $O(1)
O_{\assubs}\left(\frac{\sigma_{n,1}}{n}(\log\log n)^{1/2}\right)$.
Indeed, by the same argument as in \cite[Theorem
3]{CsorgoRevesz1978},
\begin{equation}\label{eq-estimates-on-y-theta}
\frac{y(1-y)}{\theta(1-\theta)}=O(1) .
\end{equation}
Thus, by (\ref{eq-estimates-on-y-theta}) and (B),
\begin{eqnarray*}
\lefteqn{\sup_{y\in
[\delta_n,1-\delta_n]}|(fQ)^{'}(\theta)\sqrt{y(1-y)}|}\\
&=&\sup_{y\in
[\delta_n,1-\delta_n]}|(fQ)^{'}(\theta)|(\theta(1-\theta))^{1/2}
\left(\frac{y(1-y)}{\theta(1-\theta)}\right)^{1/2}=O(1)
\end{eqnarray*}
The second order term, in view of (A), is treated in a similar way.

Consequently, by the above calculations, (\ref{step.5}) and
(\ref{conjecture}),
\begin{equation}
\sup_{y\in
[\delta_n,1-\delta_n]}|u_n(y)-\alpha_n(y)|=O_{\assubs}(a_n(\log\log
n)^{1/2}).
\end{equation}
Further,
\begin{equation}\label{eq-sup-small-uniform.quantile}
\sup_{y\in (0,\delta_n]}|u_n(y)|=O_{\assubs}(a_n (\log\log n)^{1/2})
\end{equation}
by the same argument as in \cite[Theorem 3]{CsorgoRevesz1978}. Also,
\begin{equation}\label{eq-12b}\sup_{y\in (0,\delta_n]}|\alpha_n(y)|=O_{\assubs}[\delta_n^{1-\mu}\ell
(n))+O_{\assubs}(d_{n,2})
\end{equation}
via the reduction principle, (A) and (B). Indeed,
$$
\tilde
V_{n,2}(\delta_n)=fQ(\delta_n)\sigma_{n,1}^{-1}\sum_{i=1}^nX_i+f^{'}(Q(\delta_n))\sigma_{n,1}^{-1}Y_{n,2}.
$$
The first part is $O_{\assubs}[\delta_n^{1-\mu}\ell (n))$ by (A).
For the second part, write
$$
f^{'}(Q(\delta_n))\frac{\sigma_{n,2}}{\sigma_{n,1}}=\frac{f^{'}(Q(\delta_n))}{f(Q(\delta_n))}(\delta_n(1-\delta_n))^{\mu}
\frac{f(Q(\delta_n))}{(\delta_n(1-\delta_n))^{\mu}}\frac{\sigma_{n,2}}{\sigma_{n,1}}=O(1)\frac{\delta_n^{1-\mu}}{\delta_n^{\mu}}\frac{\sigma_{n,2}}{\sigma_{n,1}}
$$
by (A) and (B). The above bound is $O(1)$ since $\mu<1/2$. Consequently, (\ref{eq-12b}) follows.

Therefore, the result of lemma follows.\koniec

From (\ref{conjecture}) with $p=2$, Lemma
\ref{lem-distance-uniform-empirical-quantile} together with the
reduction principle (\ref{approx-unif-empirical}) we conclude:
\begin{cor}\label{lem-reduction-quantile.uniform}
Under the conditions of {\rm Theorem}
{\rm\ref{thm-conv-uniform-BK}}, as $n\to\infty$,
$$\sup_{y\in
(0,1)}|u_n(y)+\sigma_{n,1}^{-1}\tilde
V_{n,2}(y)|=O_{\assubs}(a_n(\log \log n)^{1/2}),
$$
$$\sup_{y\in
(0,1)}|u_n(y)+\sigma_{n,1}^{-1}f(Q(y))\sum_{i=1}^nX_i|=O_{\assubs}(a_n(\log \log n)^{1/2}(\log n)^{1/2})
$$
and
$$
\sup_{y\in (0,1)}|u_n(y)|=O_{\assubs}((\log\log n)^{1/2}).
$$
\end{cor}
\subsection{Proof of Theorem \ref{thm-conv-uniform-BK}}
Let $\psi(y)=(y(1-y))^{\mu}$, $\mu$ from (A). Via
(\ref{approx-unif-empirical}) and (\ref{eq-replace}) we obtain
\begin{eqnarray}
\lefteqn{\sup_{y\in
[\delta_n,1-\delta_n]}\psi(y)\left|\tilde R_n(y)-\sigma_{n,1}^{-1}n^{-1}f^{'}(Q(y))\left(\sum_{i=1}^nX_i\right)^2\right|}\nonumber\\
& \le & C\sup_{y\in
(0,1)}\left|\alpha_n(y)-u_n(y)+\sigma_{n,1}^{-1}(\tilde
V_{n,2}(y)-\tilde V_{n,2}(U_n(y)))\right|\nonumber\\
&&+\sigma_{n,1}^{-1}\sup_{y\in
[\delta_n,1-\delta_n]}\psi(y)\left|(\tilde V_{n,2}(y)-\tilde
V_{n,2}(U_n(y)))+n^{-1}f^{'}(Q(y))\left(\sum_{i=1}^nX_i\right)^2\right|\nonumber\\
&=: & O_{\assubs}(d_{n,2})+I_2.\label{estimation}
\end{eqnarray}
Then
\begin{eqnarray*}
\lefteqn{\sup_{y\in [\delta_n,1-\delta_n]}\psi(y)\left|(\tilde
V_{n,2}(y)-\tilde
V_{n,2}(U_n(y)))+n^{-1}f^{'}(Q(y))\left(\sum_{i=1}^nX_i\right)^2\right|}\\
& \le  & \sup_{y\in
[\delta_n,1-\delta_n]}\psi(y)\left|\left\{f(Q(y))-f(Q(U_n(y)))\right\}\sum_{i=1}^nX_i+n^{-1}f^{'}(Q(y))\left(\sum_{i=1}^nX_i\right)^2\right|\\
&&+ \sup_{y\in
[\delta_n,1-\delta_n]}\psi(y)\left|\left\{f^{'}(Q(U_n(y)))-f^{'}(Q(y))\right\}\right||Y_{n,2}|\\
& \le  & \sup_{y\in
[\delta_n,1-\delta_n]}\psi(y)\left|(fQ)^{'}(y)(y-U_n(y))\sum_{i=1}^nX_i+n^{-1}f^{'}(Q(y))\left(\sum_{i=1}^nX_i\right)^2\right|\\
&&+\sup_{y\in
[\delta_n,1-\delta_n]}\psi(y)\frac{1}{2}|(fQ)^{''}(\theta)|(y-U_n(y))^2\left|\sum_{i=1}^nX_i\right|\\
&&+\sup_{y\in
[\delta_n,1-\delta_n]}\psi(y)|(f^{'}Q)^{'}(\theta)||y-U_n(y)||Y_{n,2}|\\
&\le & \frac{\sigma_{n,1}}{n}\sup_{y\in
[\delta_n,1-\delta_n]}\psi(y)\left|(fQ)^{'}(y)u_n(y)+\sigma_{n,1}^{-1}f^{'}(Q(y))\sum_{i=1}^nX_i\right|\left|\sum_{i=1}^nX_i\right|\\
&& +\frac{1}{2}\sup_{y\in
[\delta_n,1-\delta_n]}\psi(y)(y(1-y))|(fQ)^{''}(\theta )|
\frac{(y-U_n(y))^2}{y(1-y)}\left|\sum_{i=1}^nX_i\right|+\\
&&\sup_{y\in [\delta_n,1-\delta_n]}\psi(y)(y(1-y))^{1/2}| (f^{'}Q
)^{'}(\theta )| \frac{|y-U_n(y)|}{(y(1-y))^{1/2}}|Y_{n,2}|
\end{eqnarray*}
with the very same $\theta$ as in Lemma
\ref{lem-distance-uniform-empirical-quantile}.

As to the second term, by the condition (C) and
(\ref{eq-estimates-on-y-theta}) we have
\begin{eqnarray*}
\lefteqn{\sup_{y\in
[\delta_n,1-\delta_n]}(y(1-y))^{1+\mu}|(fQ)^{''}(\theta )|}\\
&=&\sup_{y\in
[\delta_n,1-\delta_n]}(\theta(1-\theta))^{1+\mu}|(fQ)^{''}(\theta
)|\left(\frac{y(1-y)}{\theta(1-\theta)}\right)^{1+\mu}=O(1).
\end{eqnarray*}
Thus, via Lemma \ref{lem-weighted.lil-unifom-quantile} and
(\ref{step.5}), the order of the second term is no greater than
$O_{\assubs}(\sigma_{n,1}^2n^{-2}\sigma_{n,1}(\log\log
n)^{3/2})=O_{\assubs}(n^{5/2-3\beta}\ell (n))$.

For the third term, via condition (A) and (\ref{eq-estimates-on-y-theta})
$$
(f^{'}Q )^{'}(\theta )(y(1-y))^{1/2+\mu}=(f^{'}Q )^{'}(\theta )(\theta(1-\theta))^{1/2+\mu}\left(\frac{y(1-y)}{\theta(1-\theta)}\right)^{1/2+\mu}=O(1).
$$
Consequently, the third term is $O_{\assubs}(\sigma_{n,1}n^{-1}\sigma_{n,2}\ell(n))=O_{\assubs}(n^{5/2-3\beta}\ell (n))$. \\

As for the first term, we bound this by
\begin{eqnarray*}
\lefteqn{\frac{\sigma_{n,1}}{n}\sup_{y\in
[\delta_n,1-\delta_n]}\psi(y)|(fQ)^{'}(y)|\left|u_n(y)+\sigma_{n,1}^{-1}\tilde
V_{n,2}(y)\right|\left|\sum_{i=1}^nX_i\right|}\\
&&+ n^{-1}\sup_{y\in
[\delta_n,1-\delta_n]}\psi(y)\left|f^{'}(Q(y))\sum_{i=1}^nX_i-(fQ)^{'}(y)\tilde
V_{n,2}(y)\right|\left|\sum_{i=1}^nX_i\right|=:I_3+I_4.
\end{eqnarray*}
From the condition (B), (\ref{step.5}) and Corollary
\ref{lem-reduction-quantile.uniform}, the term $I_3$ is
$$
O_{\assubs}(\sigma_{n,1}n^{-1}a_n(\log\log
n)^{1/2}\sigma_{n,1}(\log\log
n)^{1/2})=O_{\assubs}(n^{5/2-3\beta}\ell(n)).
$$
Noting that $(fQ)^{'}(y)f(Q(y))=f^{'}(Q(y))$, the term $I_4$ equals
$$
n^{-1}\sup_{y\in
[\delta_n,1-\delta_n]}\psi(y)|(fQ)^{'}(y)f^{'}(Q(y))||Y_{n,2}|\left|\sum_{i=1}^nX_i\right|=O_{\assubs}(n^{5/2-3\beta}\ell(n))
$$
since $\psi(y)(fQ)^{'}(y)f^{'}(Q(y))=O(1)$.

Thus, the term $I_2$ in (\ref{estimation}) is
$O_{\assubs}(\sigma_{n,1}^{-1}n^{5/2-3\beta}\ell(n))$. Consequently,
\begin{eqnarray*}
\lefteqn{\sup_{y\in
[\delta_n,1-\delta_n]}\psi(y)\left|\tilde R_n(y)-\sigma_{n,1}^{-1}n^{-1}f^{'}(Q(y))\left(\sum_{i=1}^nX_i\right)^2\right|}\\
&= &
O_{\assubs}(d_{n,2})+O_{\assubs}(\sigma_{n,1}^{-1}n^{5/2-3\beta}\ell (n))=O_{\assubs}(d_{n,2}).
\end{eqnarray*}
Therefore,
$$
\sup_{y\in [\delta_n,1-\delta_n]}\left|n\sigma_{n,1}^{-1}\tilde
R_n(y)-\sigma_{n,1}^{-2}f^{'}(Q(y))\left(\sum_{i=1}^nX_i\right)^2\right|=O_{\assubs}(n\sigma_{n,1}^{-1}d_{n,2}\delta_n^{-\mu})=o_{\assubs}(1)
$$
since $0<\mu<1/2$.\koniec
\subsection{Proof of Theorem \ref{thm-conv-uniform-verwaat}}
We have for $t<1/2$,
\begin{eqnarray*}
\lefteqn{2\sigma_{n,1}^{-1}n\int_0^t\tilde
R_n(y)dy=2\sigma_{n,1}^{-1}n\int_{(0,t)\cap
[\delta_n,1-\delta_n]}\tilde
R_n(y)dy}\\&&+O\left(\sigma_{n,1}^{-1}n\int_0^{\delta_n}|u_n(y)|dy\right)+O\left(\sigma_{n,1}^{-1}n\int_0^{\delta_n}|\alpha_n(y)|dy\right).
\end{eqnarray*} The second integral is at most of the order
$$
O_{\assubs}\left(\sigma_{n,1}^{-1}n\delta_n\sup_{y\in
(0,\delta_n]}|u_n(y)|\right)=o_{\assubs}(1)
$$
by (\ref{eq-sup-small-uniform.quantile}). The same holds for the
third one. A similar reasoning applies for $t>1/2$. Thus, the result
follows from Corollary \ref{thm-conv-uniform-BK-corollary}.\koniec
\subsection{Proof of Theorem \ref{thm-conv-uniform-verwaat.error}}
As in \cite{CSW2006}, let
$$
A_n(t)=2\sigma_{n,1}^{-1}n\int_{U_n(t)}^t(\alpha_n(y)-\alpha_n(t))dy.
$$
Then, $\tilde W_n(t)=A_n(t)-\tilde R_n^2(t)$ (cf. (3.7) in
\cite{CaskiCsorgoFoldesShiZitikise}). Hence, via Theorem
\ref{thm-conv-uniform-BK} and (\ref{step.5}),
\begin{equation}\label{eq-1}
\sup_{t\in [\delta_n,1-\delta_n]}|A_n(t)-\tilde
W_n(t)|=O_{\assubs}(n^{- (2\beta-1)}\ell(n)).
\end{equation}
Via the reduction principle and the second part of Corollary
\ref{lem-reduction-quantile.uniform},
\begin{eqnarray}
\lefteqn{\sup_{t\in (0,1)}|A_n(t)+B_n(t)|=:\sup_{t\in
(0,1)}\left|A_n(t)+2\sigma_{n,1}^{-2}n\int_{U_n(t)}^t(\tilde
V_{n,2}(y)-\tilde V_{n,2}(t))dy\right|}\label{eq-2}\\
& \le  & 4\sigma_{n,1}^{-1}n\sup_{y\in (0,1)}|y-U_n(y)|\sup_{y\in
(0,1)}|\alpha_n(y)+\sigma_{n,1}^{-1}\tilde
V_{n,2}(y)|=O_{\assubs}(d_{n,2}(\log\log n)^{1/2}).\nonumber
\end{eqnarray}
Let $C(t)=\int_0^tf(Q(y))dy$, $D(t)=\int_0^tf^{'}(Q(y))dy$ . Then
\begin{eqnarray*}
\lefteqn{B_n(t)=2\sigma_{n,1}^{-2}n\left(\sum_{i=1}^nX_i\int_{U_n(t)}^t(f(Q(y))-f(Q(t)))dy-Y_{n,2}\int_{U_n(t)}^t(f^{'}(Q(y))-f^{'}(Q(t)))dy\right)}\\
&=& 2\sigma_{n,1}^{-2}n\sum_{i=1}^nX_i\left\{C(t)-C(U_n(t))-(t-U_n(t))f(Q(t))\right\} -\\
&& 2\sigma_{n,1}^{-2}nY_{n,2}\left(D(t)-D(U_n(t))-(t-U_n(t))f^{'}(Q(t))\right)\\
&= & 2\sigma_{n,1}^{-2}n\frac{(fQ)^{'}(t)}{2}(t-U_n(t))^2\sum_{i=1}^nX_i+2\sigma_{n,1}^{-2}n\frac{(fQ)^{''}(\theta)}{6}(t-U_n(t))^3\sum_{i=1}^nX_i-\\
&& 2\sigma_{n,1}^{-2}n\frac{(f^{'}Q)^{'}(t)}{2}(t-U_n(t))^2Y_{n,2}-2\sigma_{n,1}^{-2}n\frac{(f^{'}Q)^{''}(\theta)}{6}(t-U_n(t))^3Y_{n,2}.
\end{eqnarray*}
where $\theta$ is from Lemma
\ref{lem-distance-uniform-empirical-quantile}. Consequently, by
(\ref{step.5}) and (\ref{conjecture}),
\begin{eqnarray}
\lefteqn{\sup_{t\in
[\delta_n,1-\delta_n]}\left|B_n(t)-\sigma_{n,1}^{-2}n
\sum_{i=1}^nX_i(fQ)^{'}(t)(t-U_n(t))^2\right|\label{eq-B_n-reduction}}\\
&= & O_{\assubs}(\sigma_{n,1}^{-1}n\ell (n))\times
\sup_{t\in
[\delta_n,1-\delta_n]}\frac{|t-U_n(t)|^3}{(t(1-t))^{3/2}}(fQ)^{''}(\theta)(t(1-t))^{3/2}+\nonumber\\
&&
O_{\assubs}(\sigma_{n,2}\sigma_{n,1}^{-2}n\ell (n))\times
\sup_{t\in
[\delta_n,1-\delta_n]}\frac{|t-U_n(t)|^2}{t(1-t)}(f^{'}Q)^{'}(t)(t(1-t))+\nonumber\\
&& O_{\assubs}(\sigma_{n,2}\sigma_{n,1}^{-2}n\ell (n))\times
\sup_{t\in
[\delta_n,1-\delta_n]}\frac{|t-U_n(t)|^3}{(t(1-t))^{3/2}}(f^{'}Q)^{''}(\theta)(t(1-t))^{3/2}
\nonumber.
\end{eqnarray}
Therefore, by (B), (C) and Lemma
\ref{lem-weighted.lil-unifom-quantile}, the bound in
(\ref{eq-B_n-reduction}) is of the order
$O_{\assubs}(n^{-(2\beta-1)}\ell (n))=O_{\assubs}(d_{n,2}\ell (n))$.
Consequently, via (\ref{eq-1}), (\ref{eq-2}),
$$
\sup_{t\in [\delta_n,1-\delta_n]}\left|\tilde
W_n(t)+\sigma_{n,1}^{-2}n\sum_{i=1}^nX_i(fQ)^{'}(t)(t-U_n(t))^2\right|=O_{\assubs}(d_{n,2}\ell
(n)).
$$
Therefore, weak convergence of $\tilde W_n(t)1_{\{t\in
[\delta_n,1-\delta_n]\}}$ follows from Corollary
\ref{lem-reduction-quantile.uniform} and central limit theorem for
partial sums $\sum_{i=1}^nX_i$. Further, by
(\ref{eq-sup-small-uniform.quantile}),
$$
\sigma_{n,1}^{-2}n^2\int_0^{\delta_n}|u_n(y)|dy=O_{\assubs}(\sigma_{n,1}^{-2}n^2\delta_n)\sup_{y\in
(0,\delta_n]}|u_n(y)|=o_{\assubs}(1)
$$
and the same holds if one replaces $u_n(y)$ with $\alpha_n(y)$.
Thus, $$\sigma_{n,1}^{-2}n^2\int_0^{\delta_n}|\tilde
R_n(y)|dy=o_{\assubs}(1).$$ Finally, $\sigma_{n,1}^{-1}n\sup_{t\in
(0,\delta_n]}\alpha_n^2(t)=O_{\assubs}(\sigma_{n,1}^{-1}nd_{n,2}\ell
(n))=o_{\assubs}(1)$. \koniec

\section*{Acknowledgment}
We would like to thank both Referees and the Associate Editor for
valuables comments which improved the paper.

\end{document}